\documentclass[journal]{IEEEtran}

\usepackage{amsmath,amssymb,amsfonts,bm}
\usepackage{booktabs,array}
\usepackage{graphicx}
\usepackage{stfloats}
\usepackage{placeins}
\usepackage[nospread,noshrink]{flushend}
\usepackage{microtype}
\usepackage{xcolor}
\usepackage{siunitx}
\usepackage{hyperref}
\usepackage[nameinlink,noabbrev]{cleveref}

\hypersetup{
 hidelinks,
 pdfauthor={Guoqing Zhao; Liang Li; Xing Li; Hao Kang},
 pdftitle={A Clustered Equivalent-Source Method with Local Response Compression and Group-Symmetry Reduction for Antenna-Array Analysis},
 pdfkeywords={Antenna arrays, equivalent sources, group symmetry, T-matrix methods}
}

\begin{document}

\title{A Clustered Equivalent-Source Method with Local Response Compression and Group-Symmetry Reduction for Antenna-Array Analysis}
\author{Guoqing Zhao, Liang Li, Xing Li, and Hao Kang
\thanks{G. Zhao, L. Li and H. Kang are with the School of Mathematical Sciences, University of Electronic Science and Technology of China (UESTC).  X. Li is with the School of Electronic Science and Engineering, UESTC. The corresponding authors are Liang Li and Xing Li (plum.liliang@gmail.com, plum\_liliang@uestc.edu.cn, xingli@uestc.edu.cn).}%
\thanks{Manuscript received October 31, 2026.}}

\markboth{IEEE ransactions on Antennas and Propagation ,~VOL.~X, NO.~X, XXX~2026}%
{Zhao \MakeLowercase{\textit{et al.}}: Equivalent-Source with Group-Symmetry Reduction}

\maketitle

\begin{abstract}
Full-wave analysis of finite antenna arrays is computationally demanding because detailed element responses must be resolved and coupled across the entire array. To reduce this cost, this paper presents a clustered equivalent-source method that combines local response compression with group-symmetry reduction. Local responses are represented in compressed equivalent-surface bases or coordinates inherited from child clusters. Symmetry constraints associated with the cluster partition and excitation reduce the coupled system, while the corresponding group transformations recover the complete response in the retained coordinates. To assess overall efficiency, we evaluate the costs of local construction, propagation, global solution, and recovery. Numerical examples involving layered and finite annular patch arrays show close agreement between the computed and reference complex far fields and shorter complete-run times than those of the corresponding symmetry-reduced direct solvers. Among the three equivalent-source implementations for the layered array, the one with the smallest coupled system has the longest runtime because local construction and propagation dominate the cost. These findings highlight the need to consider local construction and coupling costs when assessing the benefit of reducing the global system dimension.

\end{abstract}

\section{Introduction}

\begin{figure}[!tbp]
\centering
\includegraphics[width=0.98\columnwidth]{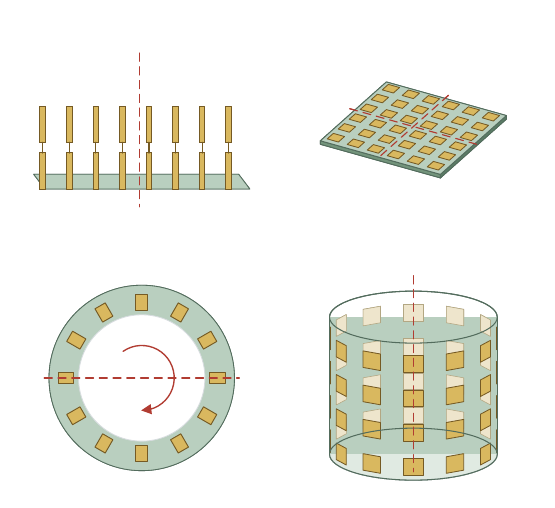}
\caption{Symmetric linear, planar, circular, and cylindrical antenna arrays.}
\label{fig:symmetric-array-geometries}
\end{figure}

Antenna arrays control the direction, gain, and spatial coverage of radiated fields by combining the responses of many elements. They support satellite communications, wireless links, and sensing, where both element geometry and excitation determine the array response \cite{BrandtMollerEtAl2024}. Mutual coupling and finite-aperture edge effects influence the response of array elements and must be represented in a finite-array model \cite{PatelEtAl2020Macromodel,Holter2002Finite}. Accurate radiation analysis therefore requires a model that connects local element responses to the interactions across the complete array. Figure~\ref{fig:symmetric-array-geometries} shows linear, planar, circular, and cylindrical arrays in the upper-left, upper-right, lower-left, and lower-right panels, respectively.

Surface integral formulations represent radiation and scattering through currents or boundary traces. The method of moments (MoM) discretizes these unknowns, commonly using Rao--Wilton--Glisson functions on conducting surfaces \cite{RaoWiltonGlisson1982}. Background Green functions incorporate the unbounded exterior \cite{SmigajEtAl2015BEM,ScroggsEtAl2017MaxwellBEM}. For a dense system with $N_{\rm dof}$ unknowns, matrix storage and direct factorization require $O(N_{\rm dof}^2)$ memory and $O(N_{\rm dof}^3)$ operations, respectively.

Fast interaction methods include the multilevel fast multipole algorithm (MLFMA) \cite{SongLuChew1997MLFMA}, its hybrids with subdomain responses \cite{ShaoEtAl2011TEPA,Garcia2008CBFMLFMA}, adaptive cross approximation \cite{Fang2020NestedACA,Jia2024ParallelACA}, precorrected fast Fourier transform (FFT) methods \cite{NieLiYuan2002PFFT}, and array decomposition \cite{BrandtMollerEtAl2024}. Finite-element and hybrid finite-element--boundary-integral methods provide alternatives for complex materials and cavity-backed antennas \cite{Monk2003Maxwell,JinVolakis1991FEBI,JinVolakis1991PatchArrays}. Local-response reduction complements these approaches by reducing the number of degrees of freedom in the coupled system.

Subdomain multilevel methods reduce the cost of planar-antenna MoM analysis \cite{SuterMosig2000SMA}. The characteristic basis function method projects mesh-level currents onto local response spaces \cite{PrakashMittra2003CBF}, with excitation-independent variants supporting reuse for different incident fields \cite{LucenteEtAl2008ExcitationIndependentCBF}. Chen et al. construct synthetic basis functions from separate feed-voltage and incident-field solution spaces using singular value decomposition \cite{Chen2012Synthetic}. Characteristic bases have also been combined with adaptive cross approximation and multilevel constructions \cite{Maaskant2008CBFACA,Laviada2009Multilevel}; model-order reduction has been applied to connected periodic arrays \cite{Shao2023MOR}. The computational cost includes both constructing these spaces and solving the reduced system.

Equivalent-source formulations transfer subdomain coupling to enclosing surfaces through equivalence-theorem domain decomposition and the equivalence principle algorithm (EPA) \cite{LiChewJiang2006EquivalenceDDM,LiChew2007EPA}. Related approaches use spherical-wave addition theorems \cite{AlianOraizi2018Multiple}, periodic Green functions \cite{Zhang2013PeriodicEPA}, and characteristic basis functions \cite{ShaoEtAl2013CBF}. Equivalent-surface macromodels reuse local responses for repeated scatterers and driven antennas \cite{PatelEtAl2017Macromodel}, with extensions to finite multilayer dielectric substrates \cite{PatelEtAl2020Macromodel}. These macromodels replace mesh-level currents on the detailed object by an equivalent-surface current description, separating the local response calculation from the array coupling problem. The enclosing surface and retained response space determine the dimensions of the local maps and the coupled system.

The $T$-matrix formalism organizes reusable single-object and multiple-scattering responses \cite{Waterman1965TMatrix,PetersonStrom1973Multiple,Mackowski1996Ensembles}. Xiang et al. combined generalized transition matrices with source-mapping vectors on reference surfaces \cite{Xiang2013GeneralizedT}. For fixed geometry, these modules and the coupling matrices can be reused as feed amplitudes and phases change. Recursive formulations construct aggregate responses for strips and patches \cite{GurelChew1990TAP,GurelChew1993TAP} and distributed obstacles \cite{DegenEtAl2023RACTMA}. The response dimension retained at each intermediate level affects the maps required to construct its parent as well as the size of the final coupled system.

Symmetry-adapted multiple-scattering methods reduce coupled $T$-matrix systems using vector transformations, object orbits, and representation spaces \cite{NecadaTorma2021}. Point-group treatments establish corresponding operator relations and block structures \cite{Schulz1999PointGroup,Kahnert2005Irreducible,Kahnert2001EBCM}; symmetry has also been incorporated in boundary-element formulations \cite{Makitalo2014Symmetry}. For clustered equivalent sources, compressed surface bases and inherited child-cluster coordinates require different local maps and construction operations. This paper examines how the cluster partition and these two response representations jointly determine the independent coefficients and the costs of local construction and global solution.

Table~\ref{tab:related-comparison} compares representative formulations, their published problem settings, response representations, symmetry treatments, and recovery operations.

\begin{table*}[!tbp]
\centering
\caption{Comparison of representative formulations and their published problem settings.}
\label{tab:related-comparison}
\footnotesize
\setlength{\tabcolsep}{3pt}
\renewcommand{\arraystretch}{1.12}
\begin{tabular}{@{}>{\raggedright\arraybackslash}p{0.14\textwidth}>{\raggedright\arraybackslash}p{0.17\textwidth}>{\raggedright\arraybackslash}p{0.23\textwidth}>{\raggedright\arraybackslash}p{0.20\textwidth}>{\raggedright\arraybackslash}p{0.20\textwidth}@{}}
\toprule
Approach & Published problem setting & Response representation and compression & Reuse and group symmetry & Orbits and response recovery \\
\midrule
Nested fast adaptive cross approximation \cite{Fang2020NestedACA} & Plane-wave scattering by perfectly conducting objects. & Full RWG current unknowns; nested low-rank compression of far-interaction matrices. & Octree hierarchy and nested factors; no point-group reduction. & Solution in the original RWG coordinates; no group-orbit recovery. \\
Recursive aggregate $T$ matrix \cite{DegenEtAl2023RACTMA} & Scattering by randomly distributed cylinders under plane-wave or external line-source illumination. & Truncated cylindrical waves; particle- and cluster-centered $T$ matrices. & Reuse across incident fields; aggregation follows addition-theorem validity criteria. & Translation operators connect particle and cluster scattering coefficients. \\
Symmetry-adapted $T$ matrix \cite{NecadaTorma2021} & Plane-wave optical scattering and lattice-mode analysis of nanoparticle arrays. & Truncated vector-spherical waves and single-particle $T$ matrices. & Point-group irreducible representations; Bloch phases for periodic systems. & Particle orbits and unitary basis transformations relate the field coefficients. \\
Present formulation & Prescribed-source radiation from finite antenna arrays. & Compressed equivalent-surface responses or inherited child-cluster coordinates. & Partition-preserving group and excitation-selected subgroup with a one-dimensional character. & Cluster-orbit representatives satisfy stabilizer constraints; group actions recover retained responses. \\
\bottomrule
\end{tabular}
\end{table*}
The cited studies in Table~\ref{tab:related-comparison} examine externally excited scattering or lattice modes, while the present benchmarks concern prescribed-source antenna-array radiation. The table compares formulation features across these problem settings. Section~\ref{sec:experiments} evaluates accuracy and computational cost against direct solvers using the same array models and excitations, without reimplementing the tabulated methods.

The main contributions are as follows.
\begin{enumerate}
\item A group-constrained clustered equivalent-source formulation is developed for symmetric antenna arrays. The coupled equations retain all inter-cluster interactions in compressed equivalent-surface bases or inherited child-cluster coordinates, with group transformations expressed in these coordinates. The partition-preserving group, cluster orbits and stabilizers, and excitation character determine the independent coefficients and the transformations that recover the full retained array response.
\item A computational complexity analysis relates the local response representation to the work of construction and global solution. It counts independent local factors throughout the hierarchy, the response columns required to construct reusable maps, and the symmetry-reduced global dimension. Basis construction, propagation, and recovery are included, so that both local construction and global dimension enter the assessment of a partition.
\end{enumerate}

\section{Group-Symmetry-Constrained Clustered Equivalent Source Method}
\label{sec:method}

The formulation uses two reductions. A local representation defines the coordinates of each cluster, and the symmetry shared by the partition and excitation determines which coefficients are independent in the global system. Table~\ref{tab:method-notation} summarizes the notation used below.

\begin{figure*}[!tbp]
\centering
\includegraphics[width=0.98\textwidth]{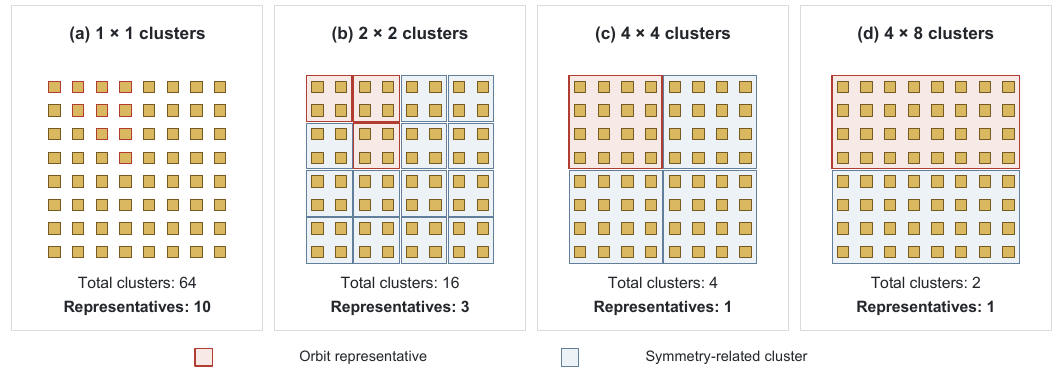}
\caption{Geometrically independent clusters for different equivalent-source partitions of an $8\times8$ array.}
\label{fig:cluster-partitions}
\end{figure*}

Figure~\ref{fig:cluster-partitions} shows how a partition changes the number of geometrically independent clusters. For the illustrated $8\times8$ array, the single-element, $2\times2$, $4\times4$, and $4\times8$ partitions contain 64, 16, four, and two clusters, with ten, three, one, and one orbit representatives marked by highlighted boundaries. The first three preserve the eight-element symmetry group $D_4$ of a square; the strip partition preserves its four-element subgroup $D_2$. The independent coefficient count also depends on retained local dimensions, stabilizers, and excitation symmetry.

\begin{table*}[!tbp]
\centering
\caption{Symbols and definitions for equivalent-source construction, symmetry reduction, and cost analysis.}
\label{tab:method-notation}
\small
\begin{tabular}{@{}p{0.2\textwidth}p{0.72\textwidth}@{}}
\toprule
Symbol & Definition \\
\midrule
$X$, $\mathcal P$, $\mathcal T$ & Finite array, complete cluster partition, and local construction hierarchy. \\
$n_B$, $m_B$ & Physical-current and equivalent-surface coordinate counts for cluster $B$. \\
$\boldsymbol i_B$, $S_B$, $\boldsymbol q_B$ & Physical-current coordinates, current-to-surface map, and equivalent-surface coefficients. \\
$L_B$, $\boldsymbol j_B$, $\boldsymbol a_B$, $T_B$, $\boldsymbol s_B$, $H_{BC}$ & Incident-to-current map, prescribed-source current response, incident coordinates, passive response, prescribed-source surface response, and propagation from $C$ to $B$. \\
$W_B$, $\Phi_{B,r}$, $\Pi_B$, $\boldsymbol c_B$, $\rho_B$, $r_B$ & Surface inner product, retained basis, left projection, retained coordinates, full response rank, and retained dimension. \\
$\mathcal D_B$, $E_{BC}$, $E_{{\rm int},B}$, $F_B$, $P_B$ & Child inverse response, external and internal reaction maps, inherited parent response, and parent-surface embedding. \\
$R_{\rm s}$, $R_{\rm a}$, $D_B(g)$ & Group actions on outgoing, incident, and retained coordinates. \\
$G_{\mathcal P}$, $\mathsf H$, $\chi$ & Partition-preserving group, excitation-selected subgroup, and its one-dimensional character. \\
$A_r$, $\boldsymbol s_r$, $\mathcal V_\chi$, $Q_\chi$, $Z_\omega$ & Compressed system and source, character space and basis, and representative-cluster stabilizer basis. \\
$d_{\rm red}$, $d_{\rm sol}$ & Exact character-space dimension and dimension of the matrix actually solved. \\
$\mathcal U(\mathcal T)$, $n_\tau$, $p_\tau$, $n_{\rm rhs}$ & Independent local factors, their dimensions and construction response-column counts, and the number of global right-hand sides. \\
\bottomrule
\end{tabular}
\end{table*}

\subsection{Coupled Equivalent-Source Formulation for a Complete Array}
\label{sec:source-construction}

Let $X$ denote the elements of a finite array. A complete partition $\mathcal P=\{B_1,\ldots,B_K\}$ consists of disjoint clusters satisfying
\begin{equation}
X=\bigcup_{k=1}^{K}B_k,\qquad B_i\cap B_j=\varnothing\quad(i\ne j).
\label{eq:complete-partition}
\end{equation}
For each cluster $B$, let $\boldsymbol i_B\in\mathbb C^{n_B}$ collect the conductor-current coordinates and, where present, equivalent-current coordinates on finite dielectric interfaces. The map $S_B\in\mathbb C^{m_B\times n_B}$ converts these physical current coordinates to equivalent-surface coefficients:
\begin{equation}
\boldsymbol q_B=S_B\boldsymbol i_B.
\label{eq:physical-to-equivalent}
\end{equation}
The equivalent currents represent the outgoing field scattered by the structures modeled within the cluster in the chosen background. Their tangential electric and magnetic components are expressed in the selected surface coordinates. The known field of the prescribed sources is treated separately.

For a linear local model, write the physical response as
\begin{equation}
\boldsymbol i_B=L_B\boldsymbol a_B+\boldsymbol j_B,
\label{eq:local-physical-response}
\end{equation}
where $\boldsymbol a_B$ contains the incident-field coordinates, $L_B$ is the passive incident-to-current map, and $\boldsymbol j_B$ is the local current response to the known fields of all prescribed sources in the absence of fields scattered by other clusters. Combining \eqref{eq:physical-to-equivalent} and \eqref{eq:local-physical-response} gives the affine equivalent-source response
\begin{equation}
\boldsymbol q_B=T_B\boldsymbol a_B+\boldsymbol s_B,\qquad
T_B=S_BL_B,\quad \boldsymbol s_B=S_B\boldsymbol j_B.
\label{eq:affine-response}
\end{equation}
The passive map $T_B$ can be reused for different prescribed sources while the cluster model is unchanged. The total radiated field is the sum of the known source field and the coupled scattered field; the direct source contribution is included once.

For distinct clusters $B$ and $C$, the propagation map $H_{BC}$ converts the outgoing response of $C$ into incident coordinates at $B$:
\begin{equation}
\boldsymbol a_B=\sum_{\substack{C\in\mathcal P\\C\ne B}}H_{BC}\boldsymbol q_C.
\label{eq:global-incident}
\end{equation}
With $\mathbf T_{\mathcal P}=\operatorname{blkdiag}(T_B)$ and the off-diagonal propagation matrix $H_{\mathcal P}$, the complete coupled system is
\begin{equation}
\begin{aligned}
A_{\mathcal P}\boldsymbol q&=\boldsymbol s_{\mathcal P},\\
A_{\mathcal P}&=I-\mathbf T_{\mathcal P}H_{\mathcal P}.
\end{aligned}
\label{eq:global-equivalent-system}
\end{equation}
The matrix $H_{\mathcal P}$ represents all interactions between distinct top-level clusters in the selected discrete surface representation. Solving \eqref{eq:global-equivalent-system} yields the coupled equivalent-source response. The local response maps, or their stored recursive factorizations, then recover physical currents or boundary traces as required.

A parent response is constructed by coupling its children internally. The parent may retain the child coefficients or express their outgoing field in a smaller parent-surface space. Both constructions preserve the separation between passive response and prescribed-source response. Section~\ref{sec:local-compression} defines these coordinate choices.

\subsection{Retained Local Response Coordinates}
\label{sec:local-compression}

Let $W_B$ define the discrete surface inner product, and let $\Phi_{B,r}$ span a retained local space of dimension $r_B$. An orthonormal representation in this inner product satisfies
\begin{equation}
\boldsymbol q_B\approx\Phi_{B,r}\boldsymbol c_B,\qquad
\Phi_{B,r}^HW_B\Phi_{B,r}=I.
\label{eq:truncated-response-basis}
\end{equation}
Substitution into \eqref{eq:global-equivalent-system} and weighted projection give
\begin{equation}
\begin{aligned}
\boldsymbol c_B&=\boldsymbol s_B^{(r)}+
\sum_{C\ne B}K_{BC}^{(r)}\boldsymbol c_C,\\
\boldsymbol s_B^{(r)}&=\Phi_{B,r}^HW_B\boldsymbol s_B,\\
K_{BC}^{(r)}&=\Phi_{B,r}^HW_BT_BH_{BC}\Phi_{C,r}.
\end{aligned}
\label{eq:compressed-coupled-system}
\end{equation}
The stacked system is denoted by $A_r\boldsymbol c=\boldsymbol s_r$, with identity diagonal blocks and off-diagonal blocks $-K_{BC}^{(r)}$. All cluster pairs remain coupled through these projected blocks. For a nonorthogonal surface basis $\Psi_B$, the corresponding left coordinate map is
\begin{equation}
\Pi_B^\Psi=(\Psi_B^HW_B\Psi_B)^{-1}\Psi_B^HW_B,
\label{eq:nonorthogonal-projection}
\end{equation}
provided that $\Psi_B$ has independent columns.

A response-adapted retained space follows from
\begin{equation}
W_B^{1/2}S_B=U_B\Sigma_BV_B^H,\qquad
\Phi_{B,r}=W_B^{-1/2}U_{B,r}.
\label{eq:weighted-response-svd}
\end{equation}
If all $\rho_B=\operatorname{rank}(S_B)$ nonzero left singular vectors are retained, the coordinates represent $\operatorname{Range}(S_B)$ exactly. Since $T_B=S_BL_B$ and $\boldsymbol s_B=S_B\boldsymbol j_B$, this gives an exact coordinate form of the chosen discrete equivalent-source system. The truncated SVD gives the best rank-$r_B$ approximation of $W_B^{1/2}S_B$ in the spectral and Frobenius norms \cite[Thms.~5.8--5.9]{TrefethenBau1997NLA}. Array-field accuracy additionally depends on the coupled response and excitation.

A parent can also retain the stacked coordinates of its children. Here, reaction coordinates are the tested incident-field entries on the right-hand side of the physical-current equations. Let $\mathcal D_B$ collect the inverse actions of the child physical systems and $\mathcal B_B$ their coordinate maps in block-diagonal form. The map $E_{{\rm int},B}$ propagates child fields into these reaction coordinates. Define
\begin{equation}
M_B=I-\mathcal B_B\mathcal D_BE_{{\rm int},B},\qquad
F_B=M_B^{-1}\mathcal B_B\mathcal D_B.
\label{eq:inherited-response}
\end{equation}
If $E_{BC}$ maps the coefficients of parent $C$ into reaction coordinates at parent $B$, the inherited system is
\begin{equation}
\boldsymbol c_B-\sum_{C\ne B}F_BE_{BC}\boldsymbol c_C
=F_B\boldsymbol p_B.
\label{eq:inherited-root}
\end{equation}
Here $\boldsymbol p_B$ stacks the prescribed-source contributions to the child right-hand sides. A separate map $P_B$ evaluates the parent-surface field as $P_B\boldsymbol c_B$. The root system couples the top-level clusters of the partition through their inherited child coordinates. The group action on these coefficients supplies the transformation used in the global reduction. Both \eqref{eq:compressed-coupled-system} and \eqref{eq:inherited-root} have the block form $A_r\boldsymbol c=\boldsymbol s_r$.

\subsection{Group-Symmetry-Constrained Global Reduction}
\label{sec:group-equivalence}

Let $G$ preserve the passive geometry, materials, element orientations, background Green function, and discrete operators. The subgroup that also preserves the complete partition is
\begin{equation}
G_{\mathcal P}=\{g\in G:g\mathcal P=\mathcal P\}.
\label{eq:partition-symmetry-group}
\end{equation}
For $g\in G_{\mathcal P}$, denote by $R_{\rm s}(g;B)$ and $R_{\rm a}(g;B)$ the transformations of outgoing surface and incident coordinates from $B$ to $gB$. Covariance of the local response and inter-cluster propagation requires
\begin{equation}
\begin{aligned}
T_{gB}R_{\rm a}(g;B)&=R_{\rm s}(g;B)T_B,\\
H_{gB,gC}R_{\rm s}(g;C)&=R_{\rm a}(g;B)H_{BC}.
\end{aligned}
\label{eq:local-covariance}
\end{equation}
These relations include the appropriate coordinate permutations and vector-component transformations.

The retained local spaces must also be closed under the group action. In the weighted-orthonormal coordinates of Section~\ref{sec:local-compression}, this condition is
\begin{equation}
R_{\rm s}(g;B)\Phi_{B,r}
=\Phi_{gB,r}D_B(g),\qquad g\in G_{\mathcal P},
\label{eq:retained-group-action}
\end{equation}
where $D_B(g)$ maps coefficients at $B$ to those at $gB$. The surface inner product must be compatible with the same transformation:
\begin{equation}
R_{\rm s}(g;B)^HW_{gB}R_{\rm s}(g;B)=W_B.
\label{eq:weight-invariance}
\end{equation}
With $\Pi_B=\Phi_{B,r}^HW_B$, closure and \eqref{eq:weight-invariance} make $D_B(g)$ unitary and give
\begin{equation}
\begin{aligned}
\Pi_{gB}R_{\rm s}(g;B)&=D_B(g)\Pi_B,\\
K_{gB,gC}^{(r)}D_C(g)&=D_B(g)K_{BC}^{(r)}.
\end{aligned}
\label{eq:projection-covariance}
\end{equation}
The second relation also uses \eqref{eq:local-covariance}. It implies that the stacked operator $A_r$ commutes with the induced coefficient action $R_r(g)$. For inherited child coordinates, the same block relation follows from covariance of the child interactions and internal solves. The induced transformations consist of unitary leaf actions and child permutations, independently of the parent-surface embedding.

Bases are mapped from orbit representatives, and singular vectors with nearly equal singular values are retained together at a truncation boundary.

The excitation selects the symmetry subspace in which the compressed system is solved. Choose a subgroup $\mathsf H\subseteq G_{\mathcal P}$ and a unit-modulus, one-dimensional character $\chi$ satisfying
\begin{equation}
R_r(h)\boldsymbol s_r=\chi(h)\boldsymbol s_r,\qquad h\in\mathsf H,
\label{eq:source-character}
\end{equation}
where $\boldsymbol s_r$ stacks the compressed source vectors. The character $\chi(h)$ is the phase factor acquired by the source under operation $h$. The invariant excitation is the special case $\chi(h)=1$ for all $h$. Define
\begin{equation}
\mathcal V_\chi=
\{\boldsymbol c:R_r(h)\boldsymbol c=\chi(h)\boldsymbol c
\ \text{for every }h\in\mathsf H\}.
\label{eq:character-space}
\end{equation}

\noindent\textbf{Proposition 1.} \textit{If $A_r$ is nonsingular, commutes with $R_r(h)$ for every $h\in\mathsf H$, and the source satisfies \eqref{eq:source-character}, the unique solution of the compressed system belongs to $\mathcal V_\chi$. If $Q_\chi$ is an orthonormal basis for this space,}
\begin{equation}
\begin{aligned}
A_\chi\boldsymbol y&=\boldsymbol s_\chi,\\
A_\chi&=Q_\chi^HA_rQ_\chi,\qquad
\boldsymbol s_\chi=Q_\chi^H\boldsymbol s_r,\\
\boldsymbol c&=Q_\chi\boldsymbol y
\end{aligned}
\label{eq:invariant-reduced-system}
\end{equation}
\textit{recovers that solution.} Appendix~\ref{app:invariant-reduction-proof} gives the proof. This restriction is exact for the retained system; local response truncation determines its approximation to the uncompressed response.

An orbit contains all clusters obtained by applying $\mathsf H$ to one cluster. For orbit-based assembly, choose one representative cluster $B_\omega$ from each orbit $\omega\in\mathcal P/\mathsf H$. Its stabilizer is
\begin{equation}
\mathsf H_{B_\omega}=\{h\in\mathsf H:hB_\omega=B_\omega\}.
\label{eq:cluster-stabilizer}
\end{equation}
Let $D_{B_\omega}(h)$ be the retained-coordinate action of a stabilizer operation, and let $Z_\omega$ span its local character space:
\begin{equation}
D_{B_\omega}(h)Z_\omega=\chi(h)Z_\omega,\qquad h\in\mathsf H_{B_\omega}.
\label{eq:stabilizer-character-space}
\end{equation}
If $hB_\omega$ is an orbit image, its coefficient vector is reconstructed as
\begin{equation}
\boldsymbol c_{hB_\omega}
=\chi(h)^{-1}D_{B_\omega}(h)Z_\omega\boldsymbol b_\omega.
\label{eq:orbit-reconstruction}
\end{equation}
The stabilizer constraint makes the lift independent of which group operation maps the representative to a given cluster. The dimension of the character-restricted system is therefore
\begin{equation}
d_{\rm red}=
\sum_{\omega\in\mathcal P/\mathsf H}\dim\operatorname{Range}(Z_\omega).
\label{eq:exact-reduced-dimension}
\end{equation}
Orthonormal stabilizer bases and an orbit-size normalization of $1/\sqrt{|\omega|}$ give an orthonormal $Q_\chi$. Each representative-target equation includes the contributions of all source clusters through their orbit maps. We use $d_{\rm sol}$ for the dimension of the factorized system; it equals $d_{\rm red}$ when all stated constraints are imposed in the solved coordinates.

After the reduced system is solved, the orbit and stabilizer maps recover all retained cluster coefficients. The surface bases evaluate the outgoing fields, and the stored local response factors recover physical currents or traces.

\subsection{Computational Complexity Analysis}
\label{sec:complexity-analysis}

The main costs are constructing the local responses, coupling the clusters, and solving the resulting system. For a construction hierarchy $\mathcal T$, we write
\begin{equation}
\begin{aligned}
\mathcal C_{\rm total}={}&\mathcal C_{\rm local}+\mathcal C_{\rm basis}
+\mathcal C_{\rm propagation}\\
&+\mathcal C_{\rm global}+\mathcal C_{\rm recovery/output}.
\end{aligned}
\label{eq:total-computational-cost}
\end{equation}
Here $\mathcal C_{\rm local}$ includes local matrix construction, prescribed-source evaluation, factorization, and response solves. $\mathcal C_{\rm basis}$ covers surface evaluation, basis extraction, and basis transformations. $\mathcal C_{\rm propagation}$ covers kernel integration, displacement templates, coordinate transformations, and coupling-block construction. $\mathcal C_{\rm global}$ covers root-matrix assembly, factorization, and solution. $\mathcal C_{\rm recovery/output}$ covers coefficient and current recovery, residuals, far fields, and output. Each operation is counted once, with stored factors and maps reused throughout the computation.

Let $\mathcal U(\mathcal T)$ denote the independent local factors. For factor $\tau$, $n_\tau$ is its matrix dimension, and $p_\tau$ is the number of construction right-hand sides solved using that factor. The local factorization and solve cost is
\begin{equation}
\mathcal C_{\rm local,solve}=O\!\left(\sum_{\tau\in\mathcal U(\mathcal T)}
\bigl(n_\tau^3+n_\tau^2p_\tau\bigr)\right).
\label{eq:local-solve-complexity}
\end{equation}
The cubic term represents factorization; the quadratic term represents solves with the stored factors. Each independent leaf, parent, or symmetry-sector matrix is factorized once, and symmetry-related copies reuse its factors. Constructing a reusable passive response includes all required symmetry sectors. The prescribed source determines the character space used in the root solve. Factor storage is $O(\sum_\tau n_\tau^2)$.

For an explicit $m_B\times n_B$ surface map $S_B$, $m_B$ and $n_B$ are the surface and physical dimensions. Forming the Gram matrix $S_BS_B^H$ costs $O(m_B^2n_B)$, and its dense eigendecomposition costs $O(m_B^3)$. Using symmetry-representative columns replaces $n_B$ by their number in the Gram-matrix cost and adds the group transformations. An inherited parent uses the child coordinates as its response representation. Surface evaluation and basis transformations complete $\mathcal C_{\rm basis}$; the response solves are included in \eqref{eq:local-solve-complexity}.

Propagation cost depends on the independent displacement templates, kernel quadrature, and response dimensions. For each independent leaf displacement, a leaf-displacement registry stores a map from source coefficients to the receiving physical right-hand side, together with the reduced passive-coupling block. Translated or symmetry-related cluster pairs reuse these entries. With $N_{\rm ker}$ templates and $r$ response columns per template, constructing the passive coupling requires $N_{\rm ker}r$ factor applications, included in $p_\tau$. Kernel integration, compression, and matrix products contribute to $\mathcal C_{\rm propagation}$. Multiplying an $r\times n$ response factor by an $n\times r$ propagated-input map costs $O(nr^2)$ per template. An inherited representation combines the stored leaf blocks in child coordinates. Factor applications used for final recovery contribute to $\mathcal C_{\rm recovery/output}$.

For a root system with dimension $d_{\rm sol}$ and $n_{\rm rhs}$ right-hand sides, dense LU factorization and solution cost
\begin{equation}
\mathcal C_{\rm global,solve}=O\!\left(d_{\rm sol}^3+
n_{\rm rhs}d_{\rm sol}^2\right),
\label{eq:global-solve-complexity}
\end{equation}
with $O(d_{\rm sol}^2)$ matrix and factor storage. The cost of inserting and accumulating the coupling blocks in the root matrix is included in $\mathcal C_{\rm global}$. Right-hand sides in the same excitation character space reuse the factorization while the passive system is unchanged. Excitation channels with independent factorizations add separate factorization and solve costs. Recovery and field evaluation depend on the numbers of reconstructed current coefficients, source channels, and observation directions.

The retained array has $\sum_Br_B$ coordinates before symmetry reduction. The excitation character and stabilizer constraints reduce this dimension according to \eqref{eq:exact-reduced-dimension}. For $K$ clusters with equal retained dimension $r$ and trivial stabilizers under $\mathsf H$,
\begin{equation}
d_{\rm red}=\frac{Kr}{|\mathsf H|}.
\label{eq:free-orbit-solve-dimension}
\end{equation}
For one character channel with $d_{\rm sol}=d_{\rm red}$, the leading cubic factorization term decreases by a factor $|\mathsf H|^3$ relative to the unreduced retained system. The quadratic storage term decreases by a factor $|\mathsf H|^2$. Equation~\eqref{eq:exact-reduced-dimension} gives the dimension for nontrivial stabilizers.

Equation~\eqref{eq:total-computational-cost} combines these contributions to describe how a cluster partition affects local construction, the global solve, and response recovery.

\section{Numerical Validation}
\label{sec:experiments}

\subsection{Experiment I: Layered Patch Array}
\label{sec:experiment-patch}
The layered patch-array study compares three equivalent-source implementations with a SCUFF-EM $D_4$ direct MoM solver. These four implementations use the same conductor mesh and prescribed excitation. Feko MoM provides an independent reference for far-field accuracy. The comparison evaluates complex far fields and recorded computational costs.

\subsubsection{Physical model and source symmetry}
The array contains $8\times8$ PEC square patches at 10 GHz. Each patch has a 10-mm side, and adjacent patch centers are spaced 15 mm apart. The patches lie at $z=0$, with centers $x,y\in15\{-3.5,-2.5,\ldots,3.5\}$ mm. A lossless dielectric with relative permittivity 2.2 occupies $-0.794<z<0$ mm above a PEC ground plane. The dielectric and ground extend infinitely in the transverse directions and enter the calculation through the layered Green function, as in the SCUFF-EM background formulation \cite{SCUFFMicrostrip}. The finite metal array is discretized by 1,976 Rao--Wilton--Glisson (RWG) functions per patch, giving 126,464 conductor unknowns.

Four prescribed, noncontact line-current sources lie beneath each patch at offsets $(1.5,0)$, $(0,1.5)$, $(-1.5,0)$, and $(0,-1.5)$ mm from its center. Each source extends in the $+z$ direction from $z=-0.65$ to $-0.15$ mm; endpoint charges are included. Their peak complex current weights are $(1,i,-1,-i)$ A under the $\exp(-i\omega t)$ convention. This defines the same quadrature drive at every patch.

\begin{figure}[!tbp]
\centering
\includegraphics[width=0.98\linewidth]{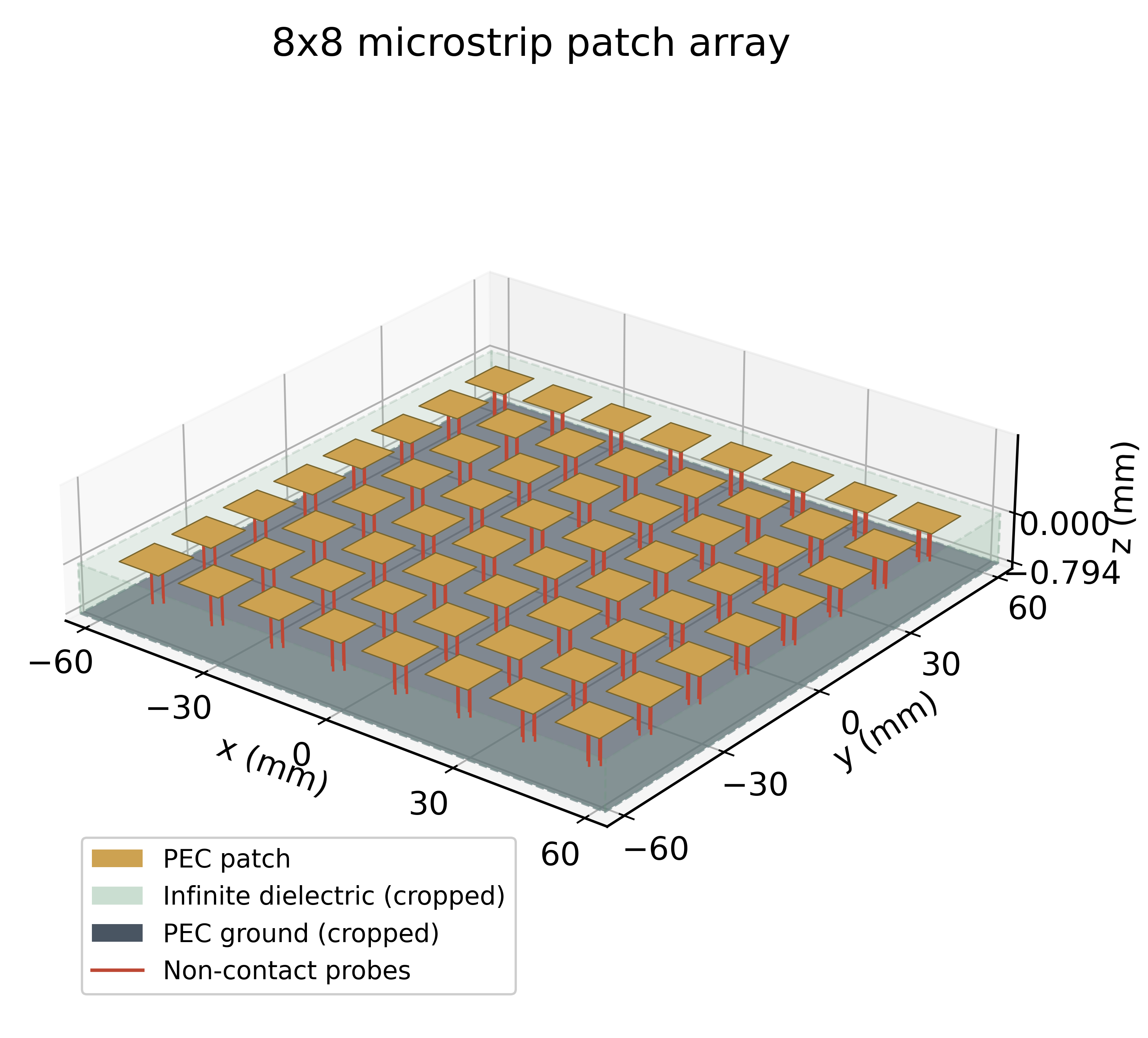}
\caption{Layered $8\times8$ patch array with prescribed noncontact line-current sources, shown schematically.}
\label{fig:patch-geometry}
\end{figure}

Figure~\ref{fig:patch-geometry} shows the finite metal array and prescribed sources. The infinite dielectric and ground are cropped for display, and vertical proportions are enlarged; coordinates are in millimeters.

The quadrature drive is formed from east-minus-west (EW) and north-minus-south (NS) source patterns. The EW response is solved in the character space of an excitation-selected four-element subgroup $\mathsf H$ of $D_4$. A quarter-turn gives the NS response, allowing the same root factorization to serve the quadrature drive. Appendix~\ref{app:patch-construction} specifies the group operations and character values.

Table~\ref{tab:patch-dimensions} distinguishes passive geometric reuse under $D_4$ from the EW auxiliary solve under $\mathsf H$. The representatives in Fig.~\ref{fig:cluster-partitions} refer to the passive geometric group. Appendix~\ref{app:patch-construction} gives the orbit counts, stabilizer constraints, and resulting root dimensions.

\begin{table*}[!tbp]
\centering
\caption{Degrees of freedom, cluster orbits, and root-system sizes for the patch-array formulations.}
\label{tab:patch-dimensions}
\small
\begin{tabular}{@{}lccccccc@{}}
\toprule
Configuration & Clusters & \shortstack{Explicit surface\\coefficients} & \shortstack{Retained degrees\\of freedom} & \shortstack{$D_4$\\orbits} & \shortstack{$\mathsf H$\\orbits} & \shortstack{Degrees of freedom\\per cluster} & \shortstack{Root\\dimension}\\
\midrule
$D_4$ direct MoM & -- & -- & -- & -- & -- & -- & 31,616\\
Feko MoM & -- & -- & -- & -- & -- & -- & 31,616\\
Single-patch surface & 64 & 76,800 & 13,056 & 10 & 16 & 204 & 3,264\\
$2\times2$ cluster & 16 & 76,800 & 13,824 & 3 & 4 & 864 & 3,456\\
$4\times4$ cluster & 4 & 19,200 & 3,556 & 1 & 1 & 889 & 889\\
\bottomrule
\end{tabular}
\end{table*}

\subsubsection{Response construction and comparison protocol}
\label{sec:patch-protocol}
The SCUFF-EM $D_4$ direct MoM solver and the three equivalent-source implementations use the same physical model, conductor mesh, prescribed sources, and layered kernel. The $D_4$ direct MoM reference uses the RWG surface-current discretization \cite{RaoWiltonGlisson1982}. Its symmetry reduction draws on group-theoretic reduction of discretized electromagnetic integral equations \cite{Cohoon1980GroupReduction}. In the present implementation, interactions are assembled using translation, $D_4$ symmetry, and reciprocity. The 31,616-dimensional EW auxiliary system is solved by in-place LU, and the complete conductor-current vector is restored.

Feko MoM is used primarily to validate far-field accuracy. Its layered-background implementation and computational accounting are specified in Appendix~\ref{app:patch-implementation}.

The retained spaces are constructed by symmetry mapping and are closed under the operations used at the root. Appendix~\ref{app:patch-construction} specifies the surface geometry, response maps, and independent local factors.

\subsubsection{Far-field accuracy and measured performance}
The upper-hemisphere complex far fields are compared with Feko using the solid-angle-weighted relative error defined in Appendix~\ref{app:patch-evaluation}. Table~\ref{tab:patch-performance} reports the errors and complete-run costs.

\begin{table*}[!tbp]
\centering
\caption{Complete-run time, peak process memory, complex-field error relative to Feko, and time ratio $t_{\rm MoM}/t$ for the five configurations. Here, $t_{\rm MoM}$ is the $D_4$ direct MoM complete-run time and $t$ is the complete-run time of each configuration.}
\label{tab:patch-performance}
\small
\begin{tabular}{@{}lcccc@{}}
\toprule
Configuration & \shortstack{Complete-run\\time (s)} & \shortstack{Peak process\\memory (GiB)} & \shortstack{Error vs\\Feko (\%)} & $t_{\rm MoM}/t$\\
\midrule
$D_4$ direct MoM & 1,720.762 & 15.729 & 0.764936 & 1.00\\
Single-patch surface & 161.466 & 0.597 & 0.766269 & 10.66\\
$2\times2$ cluster & 146.479 & 2.014 & 0.766477 & 11.75\\
$4\times4$ cluster & 401.343 & 2.370 & 0.761471 & 4.29\\
Feko MoM & 7,248.001 & 5.161 & -- & 0.24\\
\bottomrule
\end{tabular}
\end{table*}

For the solved EW auxiliary channel, the normalized full-system residuals defined in Appendix~\ref{app:patch-construction} are $1.08\times10^{-14}$, $5.40\times10^{-14}$, and $7.16\times10^{-14}$ for the single-patch, $2\times2$, and $4\times4$ configurations. These residuals measure the algebraic solution of the retained systems; the far-field error evaluates the resulting radiation against the Feko reference.

\begin{figure*}[!tbp]
\centering
\includegraphics[width=0.98\textwidth]{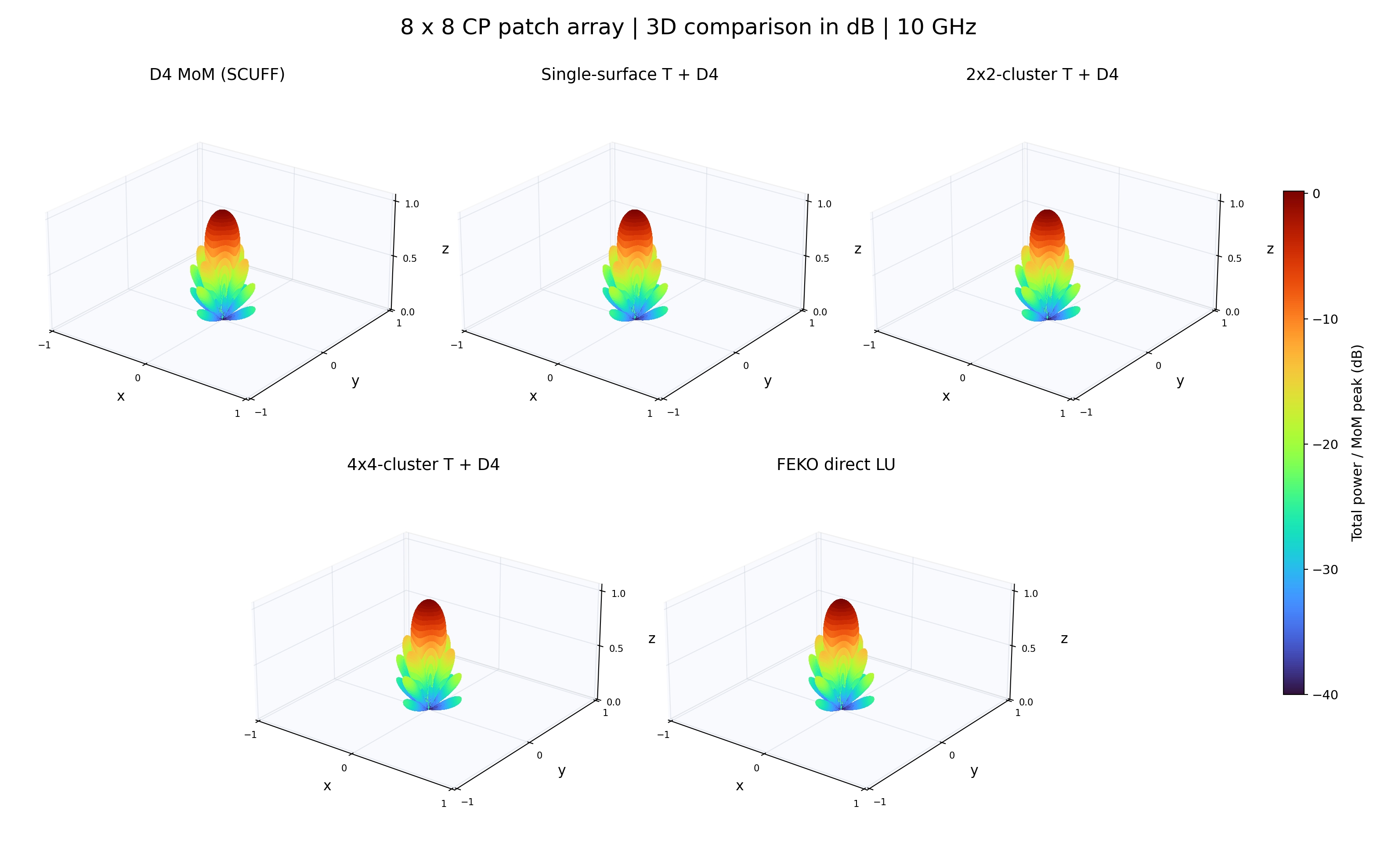}
\caption{Three-dimensional power patterns from $D_4$ direct MoM (SCUFF), three equivalent-source formulations, and Feko. Radial distance and color show total power relative to the $D_4$ direct MoM peak, in linear units and decibels, respectively.}
\label{fig:patch-patterns}
\end{figure*}

\begin{figure*}[!tbp]
\centering
\includegraphics[width=0.98\textwidth]{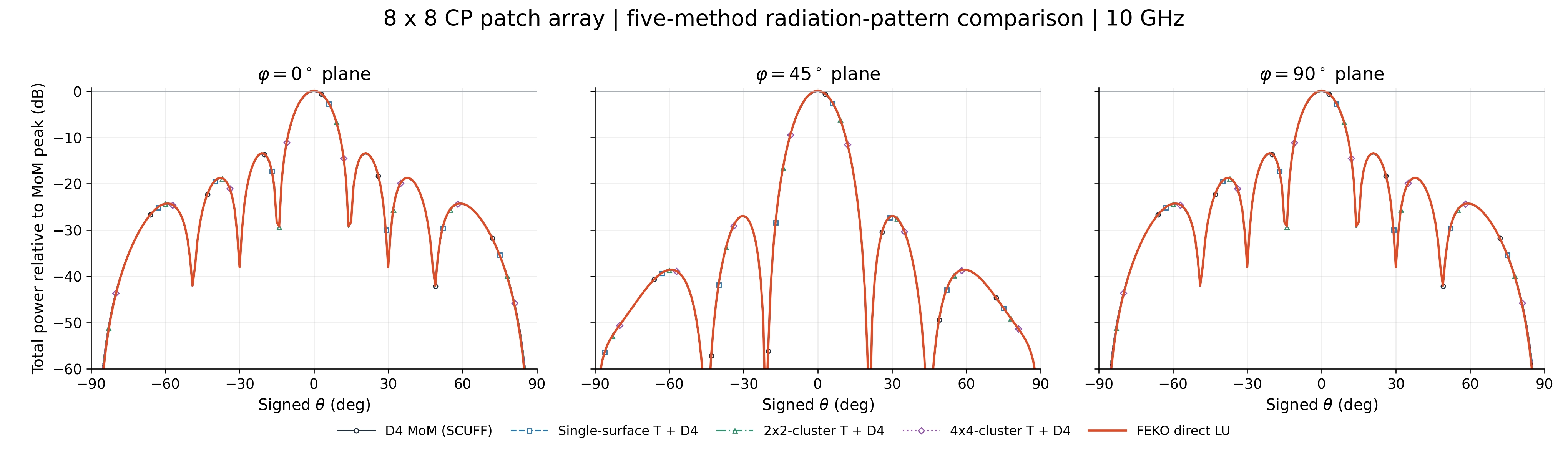}
\caption{Far-field power cuts at $\phi=0^\circ$, $45^\circ$, and $90^\circ$ for the five configurations, normalized by the $D_4$ direct MoM peak.}
\label{fig:patch-cuts}
\end{figure*}

Figure~\ref{fig:patch-patterns} compares the five three-dimensional power patterns, and Fig.~\ref{fig:patch-cuts} shows the cuts at $\phi=0^\circ$, $45^\circ$, and $90^\circ$. The SCUFF-EM $D_4$ direct MoM solution is labeled D4 MoM (SCUFF) in both figures. The five curves nearly coincide on this scale; Table~\ref{tab:patch-performance} quantifies their complex-field differences.

All four SCUFF-EM implementations yield relative complex-field errors below $0.767\%$ against Feko, with the $4\times4$ formulation giving the smallest discrepancy.

The inherited $2\times2$ implementation has the shortest recorded complete-run time, $9.3\%$ below the single-patch implementation. Source/RHS evaluation and propagation take 33.644 and 12.047 s less, respectively, offsetting 30.704 s of additional work in the other stages. 

The $4\times4$ root has $72.8\%$ fewer unknowns than the single-patch root, yet its complete run takes 2.49 times as long. Parent-field mapping and basis formation (116.588 s) and direct parent propagation (175.021 s) account for $72.7\%$ of the run. Root assembly, source response, factorization, recovery, and residual evaluation together take 5.540 s, or $1.38\%$. Parent-response construction therefore dominates the cost. Table~\ref{tab:patch-stage-times} groups disjoint driver intervals whose operations may span several terms in \eqref{eq:total-computational-cost}. Other overhead is the difference between the complete run and the listed intervals; values are rounded independently.

\begin{table*}[!tbp]
\centering
\caption{Recorded computational stage times for the three equivalent-source formulations (s).}
\label{tab:patch-stage-times}
\small
\begin{tabular}{@{}p{0.59\textwidth}ccc@{}}
\toprule
Stage & \shortstack{Single-patch\\surface} & $2\times2$ cluster & $4\times4$ cluster\\
\midrule
Local physical maps, caches, and leaf LU & 31.171 & 47.260 & 30.140\\
Leaf response basis & 1.143 & 1.228 & 0.992\\
Parent mapping and surface basis/embedding & -- & 4.650 & 116.588\\
Internal coupling and symmetry-sector construction & -- & 1.212 & 15.853\\
Propagation and passive coupling maps & 62.179 & 50.132 & 197.573\\
Root symmetry plan & 0.700 & 2.965 & 6.562\\
Prescribed-source fields and right-hand sides & 53.602 & 19.958 & 19.395\\
Root assembly, response, solve, recovery, and residual & 6.769 & 10.134 & 5.540\\
Quadrature-drive synthesis and far-field output & 4.163 & 4.642 & 2.882\\
Other mesh and output overhead & 1.740 & 4.297 & 5.818\\
\midrule
Complete-run time & 161.466 & 146.479 & 401.343\\
\bottomrule
\end{tabular}
\end{table*}

\FloatBarrier
\subsection{Experiment II: Finite Annular Patch Array}
\label{sec:annular-array}

The finite annular-array study compares a $C_{12}$ direct MoM solver with a $C_{12}$ surface $T$-matrix formulation. The comparison evaluates complex far-field accuracy and recorded computational cost.
\begin{figure*}[!tbp]
\centering
\includegraphics[width=0.98\textwidth]{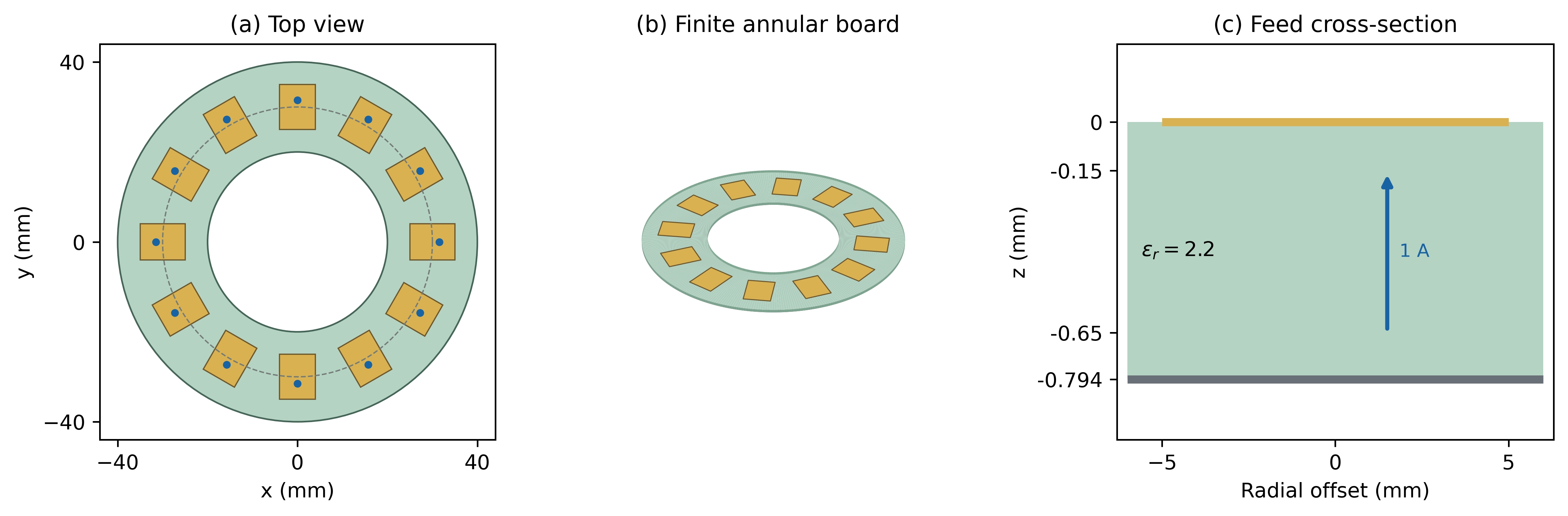}
\caption{Finite annular patch array, dielectric substrate, ground plane, and prescribed sources. The source cross section is schematic.}
\label{fig:annular-geometry}
\end{figure*}

\begin{table*}[!tbp]
\centering
\caption{Root-system degrees of freedom, complete-run time, peak process memory, and complex far-field error relative to Feko for the annular array.}
\label{tab:annular-performance}
\small
\begin{tabular}{@{}lcccc@{}}
\toprule
Configuration & \shortstack{Root degrees\\of freedom} & \shortstack{Complete-run\\time (s)} & \shortstack{Peak process\\memory (GiB)} & \shortstack{Complex-field\\error (\%)}\\
\midrule
Feko MoM & 45,294 & 10,165.90 & 1.436 & --\\
$C_{12}$ direct MoM & 15,156 & 1,404.12 & 6.947 & 0.7234\\
$C_{12}$ surface $T$ & 735 & 973.76 & 5.386 & 0.7303\\
\bottomrule
\end{tabular}
\end{table*}

\begin{figure*}[!b]
\centering
\includegraphics[width=0.98\textwidth]{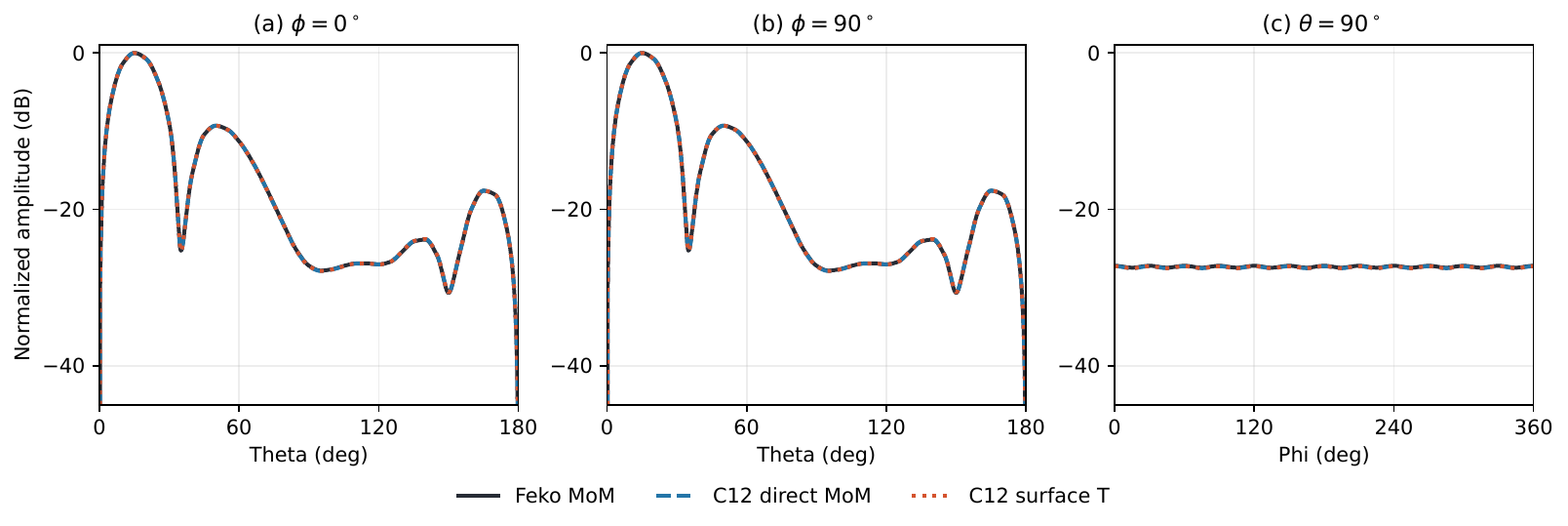}
\caption{Normalized far-field cuts at $\phi=0^\circ$, $\phi=90^\circ$, and $\theta=90^\circ$ for the annular array computed by Feko MoM, $C_{12}$ direct MoM, and $C_{12}$ surface $T$. Each pattern is normalized to its own maximum.}
\label{fig:annular-cuts}
\end{figure*}

\subsubsection{Physical model and source symmetry}
A finite annular array of twelve microstrip patches is analyzed at 10 GHz. The patch centers lie on a circle of radius 30 mm; each patch has a radial length of 10 mm and a tangential width of 8 mm. A continuous dielectric annulus with inner and outer radii of 20 and 40 mm, thickness 0.794 mm, and relative permittivity 2.2 lies above a finite PEC annular ground plane. Each patch is excited by a prescribed noncontact line current directed along $+z$, with a peak current of 1 A. The source is displaced radially outward by 1.5 mm from the patch center and extends from $z=-0.65$ to $-0.15$ mm. The sources are in phase and belong to the $m=0$ channel of the cyclic symmetry group $C_{12}$. The physical model has 181,872 unknowns. Figure~\ref{fig:annular-geometry} shows the geometry and source configuration.

\subsubsection{Response construction and comparison protocol}
The array is solved using Feko MoM with two magnetic symmetry planes, a $C_{12}$ direct MoM, and a $C_{12}$ surface $T$-matrix formulation. Feko is used primarily as a reference for assessing far-field accuracy. Its distributed out-of-core LU factorization stores part of the matrix data on disk, reducing peak process memory. The surface formulation compresses 12,288 regional surface coefficients into 735 retained modes for the symmetry-reduced root system. Appendix~\ref{app:annular-construction} specifies the response construction, far-field comparison, and computational accounting.

\subsubsection{Far-field accuracy and measured performance}
\begin{figure*}[!tbp]
\centering
\includegraphics[width=0.98\textwidth]{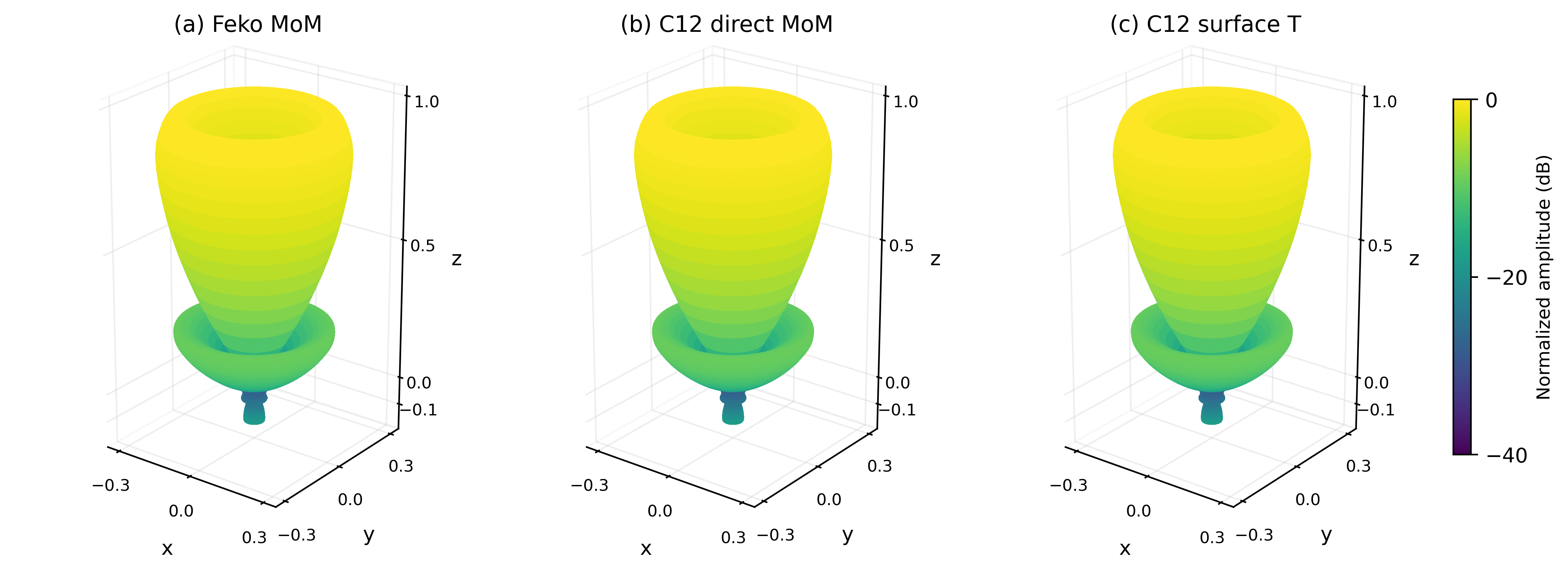}
\caption{Normalized three-dimensional radiation patterns for the annular array computed by Feko MoM, $C_{12}$ direct MoM, and $C_{12}$ surface $T$. Radial distance represents normalized field amplitude and color represents normalized amplitude in decibels.}
\label{fig:annular-patterns}
\end{figure*}

The full-sphere complex far fields are compared with the Feko reference using the solid-angle-weighted relative error defined in Appendix~\ref{app:annular-construction}. Figures~\ref{fig:annular-cuts} and \ref{fig:annular-patterns} show close agreement among the computed radiation patterns. Relative to Feko, the solid-angle-weighted complex-field errors are 0.7234\% for $C_{12}$ direct MoM and 0.7303\% for the surface $T$ formulation. Table~\ref{tab:annular-performance} compares the root-system dimensions and recorded costs. The surface $T$ run is shorter than both the $C_{12}$ direct MoM and Feko runs; relative to $C_{12}$ direct MoM, it reduces complete-run time by 30.65\% and peak process memory by 22.47\%.

\section{Conclusion}
\label{sec:conclusions}

A clustered equivalent-source method has been presented that combines local response compression with group-symmetry reduction. By separating passive scattering responses from prescribed-source responses and expressing group actions in compressed surface or inherited child-cluster coordinates, the method constructs reduced coupled systems and recovers the full retained response. The cost analysis accounts for local factorization, basis construction, intercluster propagation, the root solve, and recovery.

For the fixed $8\times8$ layered patch array, the three equivalent-source implementations yield solid-angle-weighted relative complex far-field errors of $0.7615\%$--$0.7665\%$ against Feko. The ratios of the $D_4$ direct MoM runtime to the implementation runtimes are 4.29--11.75. All use less process memory than the $D_4$ direct MoM reference. The single-patch implementation has the lowest memory use. Relative to the single-patch implementation, the $2\times2$ implementation reduces runtime by 9.3\% while using 3.37 times as much memory. Although the smallest root system reduces the global solve cost, its implementation has the longest total runtime among the three equivalent-source implementations because parent-cluster basis construction and intercluster propagation dominate. For the finite annular array, the surface $T$-matrix formulation also closely reproduces the reference far field while reducing the root-system size, complete-run time, and peak process memory relative to the $C_{12}$ direct MoM solver. Overall computational efficiency depends on the costs of local construction and propagation as well as the global system dimension.

\section*{Generative AI usage declaration}
This paper used OpenAI Codex to generate the code for the numerical implementations in Experiments~I and~II (Section~\ref{sec:experiments}). Codex was also used to plot Figs.~\ref{fig:patch-geometry}, \ref{fig:patch-patterns}, \ref{fig:patch-cuts}, \ref{fig:annular-geometry}, \ref{fig:annular-cuts}, and~\ref{fig:annular-patterns} from numerical data supplied by the authors and to polish the language describing the results in Section~\ref{sec:experiments}. The authors reviewed the experimental protocols and details and take responsibility for the code, figures, results, and final manuscript.

\section*{Acknowledgment}
This work is supported by NSFC (Grant No. 12571399 and 12101511) and Fundamental Research Funds for the Central Universities (Grant No. JBK22YJ37).

\appendices

\section{Proof of the Group-Reduction Proposition}
\label{app:invariant-reduction-proof}

For $h\in\mathsf H$, nonsingularity and equivariance imply that $A_r^{-1}$ commutes with $R_r(h)$. Hence
\[
\begin{aligned}
R_r(h)\boldsymbol c
&=A_r^{-1}R_r(h)\boldsymbol s_r\\
&=\chi(h)A_r^{-1}\boldsymbol s_r
=\chi(h)\boldsymbol c.
\end{aligned}
\]
Thus $\boldsymbol c\in\mathcal V_\chi$ and $\boldsymbol c=Q_\chi\boldsymbol y$. Equivariance makes $\mathcal V_\chi$ invariant under $A_r$, giving $A_rQ_\chi=Q_\chi(Q_\chi^HA_rQ_\chi)$. Multiplication of the compressed equation by $Q_\chi^H$ yields \eqref{eq:invariant-reduced-system}. The restriction is nonsingular on this invariant space, so the reduced solution is unique and its lift solves the compressed system.

For an orbit image $hB_\omega$, the component equation $[R_r(h)\boldsymbol c]_{hB_\omega}=\chi(h)\boldsymbol c_{hB_\omega}$ gives \eqref{eq:orbit-reconstruction}. If two operations reach the same cluster, they differ by a stabilizer operation. Its character constraint on $Z_\omega$ makes the lifted coefficient identical for either choice. Independent choices at the orbit representatives then give \eqref{eq:exact-reduced-dimension}.\hfill$\square$

\section{Layered Patch-Array Response Construction and Evaluation}
\label{app:patch-construction}

At a patch, the physical-current and equivalent-surface maps give $\boldsymbol q_\ell=S_\ell\boldsymbol i_\ell$. With patch matrix $A_\ell$ and incident-trace-to-reaction map $R_\ell$, the local factors are $T_\ell=S_\ell A_\ell^{-1}R_\ell$ and $\boldsymbol s_\ell=S_\ell A_\ell^{-1}\boldsymbol p_\ell$. The layered Green function includes the dielectric and PEC ground; the enclosing ground-backed box therefore uses explicit equivalent currents on its other five faces.

\subsection{Source symmetry and root-system dimensions}

Write the source as $\boldsymbol p_{\rm CP}=\boldsymbol p_{\rm EW}+i\boldsymbol p_{\rm NS}$, where EW and NS denote the east-minus-west and north-minus-south patterns. Let $r$ be the global $90^\circ$ rotation, $\sigma_x:(x,y)\mapsto(-x,y)$, and $\sigma_y:(x,y)\mapsto(x,-y)$. The EW source belongs to the character space of $\mathsf H=\{I,r^2,\sigma_x,\sigma_y\}$ with values $(1,-1,-1,1)$. A quarter-turn gives the NS response, allowing the same root factorization to serve the quadrature drive. 

Table~\ref{tab:patch-dimensions} distinguishes passive geometric reuse from the auxiliary solve. The $D_4$ cluster-orbit counts are ten, three, and one; the corresponding $\mathsf H$ counts are 16, four, and one. All cluster stabilizers under $\mathsf H$ are trivial in these three partitions, giving $16\times204=3264$, $4\times864=3456$, and $1\times889=889$ root degrees of freedom. The representatives in Fig.~\ref{fig:cluster-partitions} refer to the passive geometric group. For a partition with nontrivial excitation-subgroup stabilizers, the root dimension follows the local character-space sum in \eqref{eq:exact-reduced-dimension}.

\subsection{Surface coordinates and response compression}

Each implementation reuses one LU factorization of the native patch matrix $A_\ell\in\mathbb C^{1976\times1976}$. On a rectangular surface panel of area $A_p$, the polynomial functions are
\begin{equation}
\psi_{ab}(u,v)=\sqrt{\frac{(2a+1)(2b+1)}{A_p}}
\mathsf P_a(u)\mathsf P_b(v),\quad a+b\le4,
\label{eq:panel-polynomials}
\end{equation}
where $u,v\in[-1,1]$ and $\mathsf P_a$ is a Legendre polynomial. They are orthonormal under the panel area integral. Each function carries two tangential electric and two tangential magnetic components. The represented trace is $(Z_0\boldsymbol n\times\boldsymbol H_{\rm scat},-\boldsymbol n\times\boldsymbol E_{\rm scat})$. Here $\boldsymbol n$ is the outward unit normal, and $Z_0$ is the free-space wave impedance, approximately $376.73~\Omega$, used to balance the component units. Integrating the fields of each physical RWG basis function against \eqref{eq:panel-polynomials} gives a column of the stored map $S_\ell\in\mathbb C^{1200\times1976}$. Thus its coordinates already include the surface inner product.

The eigendecomposition of $S_\ell S_\ell^H$ gives the left singular space. The retained singular values satisfy $\sigma_j\ge0.003\sigma_1$, with the cutoff extended across a neighboring singular-value cluster whose relative separation is at most $10^{-7}$. With $U_\ell^HU_\ell=I$, define $B_\ell=U_\ell^HS_\ell$. The leaf response is
\begin{equation}
\begin{aligned}
\boldsymbol c_\ell&=B_\ell A_\ell^{-1}\boldsymbol p_\ell+
\sum_{j\ne\ell}B_\ell A_\ell^{-1}R_\ell
C_{\ell j}U_j\boldsymbol c_j,\\
\boldsymbol q_\ell&\approx U_\ell\boldsymbol c_\ell,
\end{aligned}
\label{eq:patch-leaf-response}
\end{equation}
where $R_\ell$ maps the incident field into physical reaction coordinates and $C_{\ell j}$ propagates the outgoing surface representation. The translated leaves share the same local basis template, with its coordinates transformed as required by the group action. The source vector $\boldsymbol p_\ell$ contains the tested fields of all prescribed line-current sources. Their known far field is added once to the conductor-scattered field.

\subsection{\texorpdfstring{Inherited $2\times2$ response}{Inherited 2-by-2 response}}

For four leaves, let $\mathcal D_B=\operatorname{blkdiag}(A_\ell^{-1})$ and $\mathcal B_B=\operatorname{blkdiag}(B_\ell)$. The inverse symbols denote applications of stored LU factors. With the internal leaf-to-reaction map $E_{{\rm int},B}$, \eqref{eq:inherited-response} gives an 864-dimensional matrix $M_B$. The parent-to-parent blocks are $F_BE_{BC}$, where $E_{BC}\in\mathbb C^{7904\times864}$ is assembled from the leaf-displacement registry. The root equations are \eqref{eq:inherited-root}; they use the four inherited leaf coefficient vectors at each parent.

The five independent $D_4$ factors of $M_B$ have dimensions 109, 107, 107, 109, and 216. Here $\mathsf E$ labels the two-dimensional irreducible representation of $D_4$, distinct from the reaction maps $E_{BC}$. Its two representation components share the same $216\times216$ multiplicity-space factor. All passive sectors enter the parent response, while the auxiliary source character is selected at the global root. Given the total reaction input $\boldsymbol f_B=\boldsymbol p_B+\sum_{C\ne B}E_{BC}\boldsymbol c_C$, physical currents are recovered by
\begin{equation}
\begin{aligned}
\boldsymbol i_B&=L_{{\rm phys},B}\boldsymbol f_B,\\
L_{{\rm phys},B}&=\mathcal D_B(I+E_{{\rm int},B}F_B),\\
\mathcal B_BL_{{\rm phys},B}&=F_B.
\end{aligned}
\label{eq:patch-internal-response}
\end{equation}
The last identity follows from $M_BF_B=\mathcal B_B\mathcal D_B$ and relates the recovered currents to the solved leaf coefficients.

Let $S_{{\rm p}\leftarrow j}$ map the physical currents of leaf $j$ to the parent-surface coefficients. Stacking the four leaf maps gives $S_{\rm p}=[S_{{\rm p}\leftarrow1}\ \cdots\ S_{{\rm p}\leftarrow4}]$. The parent-surface evaluation map is constructed separately:
\begin{equation}
P_j=S_{{\rm p}\leftarrow j}B_j^+,\qquad
P_B=[P_1\ P_2\ P_3\ P_4]\in\mathbb C^{4800\times864}.
\label{eq:patch-2x2-parent-response}
\end{equation}
Here $j$ indexes the four leaves, and $B_j^+$ selects the minimum-norm current associated with a leaf coefficient for forward field evaluation. The map $P_B\boldsymbol c_B$ is used in parent-surface far-field evaluation; root coupling and physical recovery use the inherited system above. For any physical current vector,
\begin{equation}
S_{\rm p}\boldsymbol i_B-P_B\mathcal B_B\boldsymbol i_B
=\sum_{j=1}^4S_{{\rm p}\leftarrow j}(I-B_j^+B_j)\boldsymbol i_j.
\label{eq:embedding-error}
\end{equation}
This is the field component omitted by the minimum-norm embedding. For each leaf, $\|S_{{\rm p}\leftarrow j}-P_jB_j\|_F/\|S_{{\rm p}\leftarrow j}\|_F$ is approximately $1.47\%$. For the restored EW solution, the corresponding relative 2-norm discrepancy is $0.349\%$, with the surface vectors concatenated over all parents. The embedding has a 2-norm condition estimate of $5.37\times10^4$; the character-restricted root matrix has a separate 1-norm condition estimate of 12.66. The solution-specific embedding discrepancy uses the EW parent-surface trace and a relative 2-norm, whereas the radiation diagnostic uses the total far field of the synthesized quadrature drive and the weighted norm in \eqref{eq:patch-farfield-error}.

\subsection{\texorpdfstring{Compressed $4\times4$ response}{Compressed 4-by-4 response}}

The $4\times4$ parent couples sixteen 204-coordinate leaves through a 3264-dimensional internal system. Its five independent $D_4$ factors have dimensions 410, 406, 406, 410, and 816, with the two components of the $\mathsf E$ representation sharing the same $816\times816$ multiplicity-space factor. All passive sectors contribute to the parent response.

Let $\boldsymbol s_{{\rm p},k}$ be a representative column of $S_{\rm p}$, $m_k$ its RWG orbit size, and $R_{\rm p}(g)$ the signed surface action. The parent response Gram matrix is formed as
\begin{equation}
S_{\rm p}S_{\rm p}^H=
\sum_{g\in D_4}R_{\rm p}(g)
\left(\sum_{k\in\mathcal J}\frac{m_k}{8}
\boldsymbol s_{{\rm p},k}\boldsymbol s_{{\rm p},k}^H\right)R_{\rm p}(g)^H.
\label{eq:parent-orbit-gram}
\end{equation}
The set $\mathcal J$ contains 3978 representative columns for the 31616 physical parent coordinates. Orbit weighting and the group sum reconstruct the full response Gram matrix without storing the full 4800-by-31616 map. The same singular-value criterion retains 889 coordinates in $U_{\rm p}\in\mathbb C^{4800\times889}$, with $B_{\rm p}=U_{\rm p}^HS_{\rm p}$. Its action satisfies $R_{\rm p}(g)U_{\rm p}=U_{\rm p}D_{\rm p}(g)$.

Let $L_{\rm phys}$ map reaction inputs to the internally coupled physical currents and let $R_{\rm p}^{\rm in}$ map parent incident traces to those reactions. The parent surface response is
\begin{equation}
T_{\rm p}=U_{\rm p}B_{\rm p}L_{\rm phys}R_{\rm p}^{\rm in}.
\label{eq:patch-parent-response}
\end{equation}
The reusable factor $Q_{\rm p}=B_{\rm p}L_{\rm phys}\in\mathbb C^{889\times31616}$ is constructed once. For a parent-displacement template, $C$ maps outgoing parent-surface coefficients to incident traces at the receiving parent. Integrating and compressing this propagation gives $W_c=R_{\rm p}^{\rm in}CU_{\rm p}\in\mathbb C^{31616\times889}$ for each of the two independent templates. The coupling block is then obtained as $Q_{\rm p}W_c$. This construction stores the local factors and response maps while avoiding a dense physical parent matrix.

\subsection{Surface geometry and solver implementation}
\label{app:patch-implementation}

The leaf boxes occupy $[-6,6]^2$ mm relative to each patch center. The $2\times2$ and $4\times4$ parent boxes occupy $[-13.5,13.5]^2$ and $[-29,29]^2$ mm relative to their centers. All extend from $z=-0.794$ to 3 mm. There are 20 active leaf panels or 80 active parent panels, with side panels split at the dielectric interface. Degree-four output polynomials give 1200 and 4800 explicit coordinates. Table~\ref{tab:patch-cost-inputs} gives the independent factors and propagation paths.

Feko uses double-precision MoM with an infinite layered substrate and a PEC ground plane. Electric and magnetic mirror symmetries reduce the 126,464 physical RWG unknowns to a 31,616-dimensional system, which is solved by direct LU factorization. The complete-run time is 7,248.001 s, including numerical construction, solution, recovery, and far-field evaluation from the accepted mesh. The synchronized aggregate peak RSS over all processes is 5.541161 GB, equivalent to approximately 5.161 GiB as reported under peak process memory in Table~\ref{tab:patch-performance}.

\begin{table*}[!tbp]
\centering
\caption{Independent local factors, parent representations, and propagation templates for the equivalent-source formulations.}
\label{tab:patch-cost-inputs}
\small
\begin{tabular}{@{}l >{\centering\arraybackslash}p{0.25\textwidth} p{0.30\textwidth} p{0.20\textwidth}@{}}
\toprule
Configuration & Independent internal factors & Parent representation & Propagation templates\\
\midrule
Single-patch surface & -- & Leaf $1200\times1976$ response map. & 35 leaf displacements.\\
$2\times2$ cluster & 109, 107, 107, 109, 216 & Inherited child space; $4800\times864$ forward embedding. & 35 leaf displacements.\\
$4\times4$ cluster & 410, 406, 406, 410, 816 & $4800\times4800$ Gram matrix from 3978 representative columns. & Nine internal leaf and two direct parent displacements.\\
\bottomrule
\end{tabular}
\end{table*}

\subsection{Far-field comparison and residual evaluation}
\label{app:patch-evaluation}

The complex far field is evaluated through $F_\theta$ and $F_\phi$ at 32,041 distinct directions: $\theta=0^\circ,\ldots,89^\circ$ and $\phi=0^\circ,\ldots,359^\circ$, with one polar sample. For a tested configuration, the relative error with respect to Feko is
\begin{equation}
\epsilon_{\rm patch}=
\left[\frac{\sum_\ell w_\ell\|\boldsymbol F^{\rm s}_\ell-
\boldsymbol F^{\rm Feko}_\ell\|_2^2}
{\sum_\ell w_\ell\|\boldsymbol F^{\rm Feko}_\ell\|_2^2}\right]^{1/2},
\label{eq:patch-farfield-error}
\end{equation}
where $\boldsymbol F=(F_\theta,F_\phi)$ and superscript $\rm s$ denotes the tested SCUFF-EM configuration. For a nonpolar sample, $w_\ell=\Delta\phi[\cos\theta_- -\cos\theta_+]$, with $\Delta\phi=\pi/180$, $\theta_-=\theta-0.5^\circ$, and $\theta_+=\theta+0.5^\circ$. The last ring at $89^\circ$ extends to $90^\circ$, and the unique pole has weight $2\pi(1-\cos0.5^\circ)$. These weights sum to $2\pi$. The complex fields are compared without amplitude or phase alignment.

For display, let $P=|F_\theta|^2+|F_\phi|^2$ and let $P_{\max}^{\rm M}$ be the peak value of the SCUFF-EM $D_4$ direct MoM power pattern over the saved angular grid. All five patterns use this common reference. In Fig.~\ref{fig:patch-patterns}, radial distance is $P/P_{\max}^{\rm M}$, and color is $10\log_{10}(P/P_{\max}^{\rm M})$; the Cartesian axes are dimensionless. Figure~\ref{fig:patch-cuts} uses the same power normalization in decibels. The three-dimensional plots use $2^\circ$ azimuthal sampling, while complex-field errors use the saved $1^\circ$ grid. Negative signed polar angles in a cut use the corresponding $\phi+180^\circ$ direction.

For the solved EW auxiliary channel, the normalized full-system residual is $\|A_r\boldsymbol c-\boldsymbol s_r\|_2/\|\boldsymbol s_r\|_2$ after the orbit coefficients are restored.

\section{Annular Patch-Array Response Construction and Evaluation}
\label{app:annular-construction}

\subsection{Cyclic reduction and surface response}

For the annular array in Section~\ref{sec:annular-array}, the dielectric is lossless, the conductors are modeled as zero-thickness PEC, and the local mesh size at PEC edges is 0.30 mm. The prescribed line-current model includes endpoint charges and uses the $\exp(-i\omega t)$ convention. Let $g$ be a rotation of $30^\circ$ about the array axis. The twelve sectors are related by $C_{12}=\{g^k:0\leq k<12\}$. The in-phase sources select the $m=0$ channel, for which $\chi(g^k)=1$. The physical system contains 181,872 unknowns, corresponding to 15,156 coordinates in one cyclic sector. The $C_{12}$ direct MoM formulation solves this sector system, and the responses in the remaining sectors follow from the cyclic rotations.

The surface $T$ formulation represents the regional response in 12,288 equivalent-surface coordinates and retains 735 modes for the symmetry-reduced root system. Its construction also uses one LU factorization of the 15,156-dimensional exact near-neighbor system. The 735-dimensional root solve and the near-neighbor factorization are separate stages. After the root solution, the retained responses are restored through the cyclic group actions and the physical response is recovered for far-field evaluation.

\subsection{Far-field comparison and computational accounting}

The complex far fields are evaluated at $\theta=0^\circ,5^\circ,\ldots,180^\circ$ and $\phi=0^\circ,5^\circ,\ldots,355^\circ$. Let $\boldsymbol E_\ell=(E_{x,\ell},E_{y,\ell},E_{z,\ell})$ denote the complex Cartesian far-field vector at a grid sample, and let $w_\ell$ be its solid-angle weight. The full-sphere relative error with respect to Feko is
\begin{equation}
\epsilon_{\rm ann}=
\left[
\frac{\sum_\ell w_\ell
\|\boldsymbol E_\ell-\boldsymbol E_\ell^{\rm Feko}\|_2^2}
{\sum_\ell w_\ell\|\boldsymbol E_\ell^{\rm Feko}\|_2^2}
\right]^{1/2},
\quad \sum_\ell w_\ell=4\pi .
\label{eq:annular-farfield-error}
\end{equation}
With $\Delta\theta=\Delta\phi=\pi/36$, the solid-angle weight at $(\theta_i,\phi_j)$ is
\begin{equation}
\begin{aligned}
w_{ij}=\Delta\phi\bigl[&
\cos\max(0,\theta_i-\Delta\theta/2)\\
&-\cos\min(\pi,\theta_i+\Delta\theta/2)
\bigr].
\end{aligned}
\label{eq:annular-solid-angle-weight}
\end{equation}
The grid includes all azimuth samples at each pole; their weights are included in the sums. The complex fields are compared without amplitude or phase alignment. For the plotted cuts, each solution is normalized to its own peak over the full angular grid. In the $\phi=0^\circ$ and $90^\circ$ cuts, the root-mean-square difference between the independently normalized surface-$T$ and Feko amplitudes in decibels is approximately $0.0121$ dB, evaluated where the normalized Feko amplitude is at least $-20$ dB.

The calculations use an Intel Core i7-12700 host. Feko 2022.2 runs on Windows with four MPI processes and five threads per process; the two $C_{12}$ implementations run under WSL with 20 threads. The complete-run time starts from the accepted mesh and includes numerical construction, solution, recovery, and far-field evaluation. CAD and mesh generation, compilation, and plotting are excluded. The reported peak process memory is the peak aggregate resident set size of the calculation's process tree.
\subsection{Feko implementation and memory accounting}
\label{app:feko-implementation}

Feko uses magnetic symmetry planes at $x=0$ and $y=0$, corresponding to the $YZ$ and $XZ$ planes, respectively. The complete model contains 181,872 physical unknowns, and the symmetry-reduced Feko system contains 45,294 unknowns. The mesh of the symmetry-reduced computational domain contains 61,512 triangles, with 38,760 on metal surfaces and 22,752 on dielectric interfaces.

The distributed out-of-core LU factorization stores part of the matrix data on disk. The configured solver-memory limit is 8 GiB in total, or 2 GiB per MPI process. The native Feko log reports a matrix disk footprint of 30.767578 GiB, approximately 30.768 GiB.

The synchronized aggregate peak resident set size (RSS) is obtained by summing the process RSS values at each sampling time and taking the maximum over time. For the full job, including the driver, this peak is 1.435699 GiB, reported as 1.436 GiB under peak process memory in Table~\ref{tab:annular-performance}. For the solver processes alone, excluding the driver, the synchronized aggregate peak RSS is 1.379044 GiB.

\bibliographystyle{IEEEtran}
\bibliography{references}

\end{document}